\documentclass[preprint,12pt]{elsarticle}

\usepackage{graphicx, xcolor}
\usepackage{amsmath, amssymb}

\usepackage{subfigure}
\usepackage{hyperref}

\newcommand{\beq}{\begin{equation}}
\newcommand{\eeq}{\end{equation}}
\newcommand{\bsub}{\begin{subequations}}
\newcommand{\esub}{\end{subequations}}
\newcommand{\bsube}{\begin{subequations} \begin{eqnarray}}
\newcommand{\esube}{\end{eqnarray} \end{subequations}}
\newcommand{\bea}{\begin{eqnarray}}
\newcommand{\eea}{\end{eqnarray}}

\newcommand{\RR}{\mathbb{R}}

\newcommand{\PP}{\mathbb{P}}
\newcommand{\ind}{\mathbf{1}}
\newcommand{\tx}{\tilde x}
\newcommand{\tr}{\tilde r}
\newcommand{\td}{\tilde d}
\newcommand{\dif}{\mathrm{d}}

\usepackage{lineno}

\begin{document}

\begin{frontmatter}


\title{Resampling with neural networks for stochastic parameterization in multiscale systems}



\author[1,2]{Daan Crommelin}
\author[1]{Wouter Edeling}

\address[1]{Centrum Wiskunde \& Informatica, Scientific Computing Group\\
  Science Park 123, 1098 XG Amsterdam, The Netherlands}
\address[2]{Korteweg-de Vries Institute for Mathematics, University of Amsterdam \\
  Science Park 105-107, 1098 XG Amsterdam, The Netherlands}

\begin{abstract}
In simulations of multiscale dynamical systems, not all relevant processes can be resolved explicitly. Taking the effect of the unresolved processes into account is important, which introduces the need for paramerizations.  We present a machine-learning method, used for the conditional resampling of observations or reference data from a fully resolved simulation. It is based on the probabilistic classification of subsets of reference data, conditioned on macroscopic variables. This method is used to formulate a parameterization that is stochastic, taking the uncertainty of the unresolved scales into account. We validate our approach on the Lorenz 96 system, using two different parameter settings which are challenging for parameterization methods.
\end{abstract}


\end{frontmatter}


\section{Introduction}

For modeling and simulation of multiscale systems, a central problem is how to represent processes with small spatial scales and/or fast timescales. Simulations that include all relevant scales in an explicit way are often computationally too expensive, posing major challenges to the use of numerical models in a variety of disciplines (e.g. physics, chemistry, climate science, biology, engineering). However, for many problems, the model output of interest involves only large-scale model variables, therefore a natural approach to make simulations less expensive is to derive reduced models for the large-scale variables only. The effect of unresolved (small-scale) degrees of freedom on the large-scale variables must be parameterized, in order to obtain a reduced model that forms a closed system.

Constructing a parameterization (or model closure) can be approached in different ways, using e.g. analytical or computational methods \cite{weinan2007heterogeneous,pavliotis2008multiscale,kevrekidis2009equation,weinan2011principles}. Here we focus on data-driven approaches, where data of (the effect of) small-scale processes is used to infer a parameterization (e.g. \cite{sarghini2003neural,krasnopolsky2006complex,crommelin2008subgrid,nimsaila2010markov,dorrestijn2013stochastic,ma2015using,verheul2016data,lei2016data,lu2017data}). These data can stem from e.g. fully resolved (high resolution) simulations on a limited space/time domain, or from physical measurements. Data-driven methods can be particularly useful when there is no clear separation between small/fast scales and large/slow scales, so that analytical or computational approaches that rely on such a scale gap do not apply.

The last few years have seen a surge of interest in data-driven approaches, due to the rapid developments in machine learning (ML). In e.g. \cite{pathak2018hybrid,brenowitz2018prognostic,rasp2018deep,maulik2019subgrid,bolton2019applications} (and many more studies), ML techniques are proposed for parameterizing unresolved processes based on data. Although the specific ML  techniques vary, in nearly all proposed methods the parameterization is effectively deterministic. However, uncertainties in subgrid-scale responses are important, and these can be accounted for by formulating parameterizations that are stochastic rather than deterministic \cite{palmer2001nonlinear,berloff2005random,berner2017stochastic,palmer2019stochastic}. The relevance of stochastic formulations for reduced models can be understood theoretically from the Mori-Zwanzig formalism for model reduction of dynamical systems \cite{chorin2009stochastic,lu2017data,lin2019data}.

Using ML methods for stochastic parameterization has hardly been explored yet. A recent exception is  \cite{christensen2020machine}, where an approach using generative adversarial networks (GANs) is proposed. Here we follow another approach to ML-based stochastic parameterization, by combining neural network-based probabilistic classification with resampling/bootstrapping, building on recent work on parameterization with resampling \cite{verheul2016data,verheul2017covariate,edeling2019towards}.

\section{Preliminaries}

In this study we consider multiscale dynamical systems, represented by a set of coupled nonlinear ordinary differential equations (ODEs) for the time-dependent variables $x(t)$ and $y(t)$:
\bsub
\label{eq:contin1}
\bea
\tfrac{d}{dt} x & = & f(x,\sigma) \label{eq:contin1a} \\
\tfrac {d}{dt} y &  = & g(x,y) \label{eq:contin1b} \\
\sigma & = & \sigma(y)
\eea
\esub
with $t \in \RR^+_0$,  $x \in \RR^N$, $y \in \RR^K$ and $\sigma: \RR^K \mapsto \RR^N$. The initial conditions are denoted $x(0) = x_0$ and $y(0)=y_0$.

This coupled system can be the result of e.g.\ spatial discretization of a partial differential equation. In (\ref{eq:contin1}), $x$ denotes ``macroscopic" variables, representing phenomena with large spatial scales and/or long timescales. The $y$ denotes ``microscopic" variables, typically with small or fast scales. Despite these shorthand names, we do not assume strong scale separation between $x$ and $y$, neither in space nor in time.
Often, $\text{dim}(x) =: N \ll K:=\text{dim}(y)$. 
The coupling from ``micro'' to ``macro'', represented by $\sigma(y)$, can be as simple as $\sigma(y)=y$, however there are many cases where coupling effectively takes place via a quantity (e.g. a discretized spatial field) with the dimension of $x$ rather than of $y$. We note that $\sigma$ is implicitly dependent on $x$, because $y$ is coupled to $x$ via equation (\ref{eq:contin1b}).

As already discussed in the introduction,
even though the interest is often only in the behavior of $x$, one needs to simulate both $x$ and $y$ as they are coupled. This can make numerical simulations very expensive. We aim to reduce the computational cost by replacing $\sigma$ by a data-inferred approximation or surrogate $\tilde \sigma^\text{surr} (x)$ in (\ref{eq:contin1a}) so that a closed, reduced system is obtained involving only $x$. In this study, we discuss a data-driven approach in which the surrogate is not only $x$-dependent but it is also  stochastic and has memory. The stochasticity and memory-dependence of the surrogates both follow from the Mori-Zwanzig theory, see e.g. \cite{chorin2009stochastic,lin2019data} for more details.

We focus here on the case where $\sigma$ enters as an additive term for $x$, i.e. $dx/dt = f(x) + \sigma$ (we note that the methodology presented here can in principle also be used when $\sigma$ enters as a multiplicative term, however we do not perform tests on such cases here). Furthermore, we consider the discrete-time version of this system, resulting from using an (explicit) numerical time integration scheme for (\ref{eq:contin1}) with time step $\Delta t$. We define $x_j := x(t_j)$ and $y_j := y(t_j)$. For simplicity we assume a constant time step so that $t_j = j \, \Delta t$ and thus $x_0$ and $y_0$ as defined here coincide with the initial conditions defined earlier. We denote the discrete-time system by
\bsub
\label{eq:discrete}
\bea
x_{j+1} & = & F(x_{j}) + r_{j} \label{eq:discrete_x} \\
y_{j+1} & = & G(x_{j},y_{j}) \\
r_j & := & r(y_j)
\eea
\esub
Clearly, $j \in \mathbb{N}$ is the time index here. The precise form of $F$ and $G$ depends on the time integration scheme used. For instance, a simple forward Euler scheme would result in $F(x_j) = x_j + \Delta t \, f(x_j)$, $G(x_j,y_j) = y_j + \Delta t \, g(x_j,y_j)$ and $r(y_j) = \Delta t \, \sigma(y_j)$. 

The model structure in (\ref{eq:discrete}) reflects a modular computational set-up that can be encountered in various applications (see e.g.\ \cite{hoekstra2014multiscale,jansson2019regional}), where different submodels (or model components) are used for macroscopic and microscopic processes. For example, a single macromodel can be coupled to multiple micromodels that each cover a different physical process, or a different part of the spatial domain. At every $\Delta t$ timestep of the macromodel for $x$, the micromodel(s) is called to propagate $y_j$ to $y_{j+1}$ given $x_j$ (possibly using multiple smaller``micro" time steps in order to cover one ``macro" time step $\Delta t$) so that $r_j$ can be updated to $r_{j+1}$.

In terms of the discrete-time system (\ref{eq:discrete}), the aim of reduced multiscale modeling with a surrogate is to construct a computationally cheap approximation for updating $r_j$ to $r_{j+1}$, so that the expensive micromodel is no longer needed to simulate $x$.

\section{A stochastic surrogate model from data \label{sec:stoch}}

We assume that we have observation data available in the form of time series of $(x_j,r_j)$ generated by the full multiscale model (\ref{eq:discrete}). We denote these observations by $(x_j^o, r_j^o)$, $j=0,1, ..., T$ (note that the data include $r$ but not $y$). The case we have in mind here is where these data come from numerical simulation of the full multiscale model, e.g. simulation on a limited spatial domain or over a limited time interval. However they can also come from physical experiments or measurements. We assume that there is no significant observational error. 

Key to our approach is that we aim to build a surrogate model for the time evolution of $r$ by sampling from the distribution of $r_{j+1}$ conditional on the past states of $x$ and $r$, i.e. sampling from the conditional distribution
\beq
r_{j+1} \,  | \, r_j, r_{j-1}, ..., x_j, x_{j-1}, ... 
\eeq
It is usually not known how to obtain this distribution in a systematic way from the model (\ref{eq:contin1}) or (\ref{eq:discrete}). Therefore we make use of the observations $(x_j^o, r_j^o)$ to build a surrogate.

We note that we do not need to have an explicit expression for the conditional distribution of $r_{j+1}$ (or an approximation of it), we merely need to be able to sample from it. This can be achieved by \textit{resampling} from the observations in an appropriate way. Because we use resampling, we do not have to assume specific structural properties of the conditional distribution of $r$ for constructing the surrogate.


Thus, we construct a ``stochastic surrogate'' for $r$ by random sampling from $r_{j+1}^o \, | \, r_j^o, r^o_{j-1}, ..., x^o_j, x^o_{j-1}, ...$. To make this practical, we assume finite memory (i.e., $r_{j+1}^o$ does not dependent on $r_{j'}^o$ if $j-j'$ is large enough), and we define the feature vectors
\bsub
\label{eq:features}
\bea
\td_j & := & (\tr_j, \tr_{j-1}, ..., \tr_{j-J}, \tx_j, \tx_{j-1}, ..., \tx_{j-J}) \\
d^o_j & := & (r^o_j, r^o_{j-1}, ..., r^o_{j-J}, x^o_j, x^o_{j-1}, ..., x^o_{j-J})
\eea
\esub
for some finite memory depth $J$.
The feature vectors take values in the feature space, which in this case has dimension $2N(J+1)$ since $\dim (\tx) = \dim (\tr) = N$. It is straightforward to change the definition of the feature vectors, for example using different memory depth (history) for $\tx$ and $\tr$ (in (\ref{eq:features}) it is $J$ for both), including only every $n$-th time step, using functions of $\tr$ or $\tx$, or leaving out either $\tr$ or $\tx$ enitrely. For ease of exposition, we stick to the definitions in (\ref{eq:features}).

We sample $\tr_{j+1}$ randomly from the set $S_{j+1}$ consisting of all observations $r_{i+1}^o$ whose associated feature vector $d^o_i$ is close to $\td_j$. 
Our reduced model then is
\bsub
\label{eq:reduced2}
\bea
\tx_{j+1} & = & F(\tx_{j}) + \tr_{j} \label{eq:reduced2a} \\
\tr_{j+1} & : & \text{random sample from }\, S_{j+1} := \left\{ \forall \,\, r_{i+1}^o \, \left| \, \td_j \,\,\, \text{close to} \,\,\, d^o_i  \, \right. \right\} \label{eq:reduced2b}
\eea
\esub

The resampling step to update $\tr_j$ is very similar in philosophy to the local bootstrap for Markov processes proposed by \cite{paparoditis2002local} and the nearest neighbor (k-NN) resampling scheme by \cite{lall1996nearest}. By the local bootstrap procedure, ``pseudo-time series" can be generated that reproduce the temporal dependence properties of observations from a stationary Markov process of order $p$, see \cite{paparoditis2002local}.

However, an important difference with the situation considered in \cite{paparoditis2002local,lall1996nearest} is that here, we resample a quantity ($r$) that is (two-way) coupled with another quantity ($x$) which is \emph{not} resampled. We update $\tr_j$ by resampling from $r_{i+1}^o \, | \, r_i^o, x^o_i, r^o_{i-1}, x^o_{i-1}, ...$ but we update $\tx_j$ by using the model (\ref{eq:reduced2a}). In \cite{paparoditis2002local,lall1996nearest} there is no such coupling to another model involved. 

Note that we do not use the expectation (sample average) of $S_{j+1}$ to update $\tr_j$ in (\ref{eq:reduced2b}). Using the expectation may be well-suited for a one-step ahead prediction of $\tilde x_{j+1}$, however it misses the inherent uncertainty of $r_{j+1}$ given its current and past states (as summarized in the feature vector). Another aspect of sampling versus averaging is that the expectation may be ``unphysical" (i.e., yield a value not consistent with the full model (\ref{eq:discrete})) whereas the individual samples were generated by the full model and are thus entirely consistent with it. 

In previous studies \cite{verheul2016data,verheul2017covariate} we used resampling (albeit formulated somewhat differently) with an implementation relying on ``binning'': the feature space was divided into non-overlapping bins (or cells), and $\td_j$ was considered close to $d^o_i$ when both were in the same bin. At any time during simulation with the reduced model, the feature vector $\tilde{d}_j$ fell within a single bin. From the set of observation feature vectors $d_i^o$ that fell in the same bin, we randomly selected one and used the associated $r_{i+1}^o$ as $\tilde{r}_{j+1}$. 

The results with this implementation were positive, but a drawback is that one quickly runs into curse of dimension problems when binning the feature space. If we use $N_{bins}$ bins in each dimension of the feature space, the total number of bins is $N_{bins}^{2N(J+1)}$. This number grows exponentially with both $N$ and $J$. 

However, instead of binning the \emph{input} features ($d_i^o$), we can also bin the \emph{output} data ($r_{i+1}^o$). The advantage is that the curse of dimensionality can be avoided here, since it involves only one variable and no time-lagged quantities are included (see section \ref{sec:stochsurr} for further discussion). The disadvantage is that, unlike in the previous approach (where the features were binned), no simple map from $\td_j$ to the output bin (from which to sample $\tr_{j+1}$) exists. However, we can learn this mapping from the data using a neural network. This is discussed in the next section.


\section{Resampling by neural network classification}
\label{sec:stochsurr}

We introduce here an approach for resampling by combining binning and probabilistic classification using a neural network. With this approach we can generate $\tr_{j+1}$ by resampling from the observation data, as in (\ref{eq:reduced2}), without being hampered by curse of dimension problems that occur if the input feature vector is of high dimension.

The basic idea is the following. As mentioned briefly in the previous section, we discretize the space of the observation $\{ r^o_j \}$ by defining a set of $M$ non-overlapping subsets $\{ B_m \}_{m=1}^M$ (referred to as ``bins'' here)  so that $r^o_j \in B_1 \cup B_2 \cup ... \cup B_M$ for all $j$. Then we train a neural net to map the feature vector $d_{j}^o$ to a probability distribution over the bins $\{ B_m \}$. This distribution corresponds to the probabilities that $r_{j+1}^o$ sits in the various bins $B_m$, given the feature vector $d_{j}^o$. Denoting the probability distribution as the $M$-dimensional vector $\rho := (\rho_1, ..., \rho_M)^T$, 
we thus want to train a neural net $\rho^\text{NN}$ such that
\beq
\rho^\text{NN}_m (d_{j}^o) \approx \PP (r_{j+1}^o \in B_m \, | \, d_{j}^o ) \, . 
\eeq
Obviously,  $\rho_m^\text{NN} (d) \in [0,1]$ for all $m$, and $\sum_{m=1}^M \rho_m^\text{NN} (d) = 1$ for any feature vector $d$.

In the reduced model (\ref{eq:reduced2}), given $\td_{j}$, we generate $\tr_{j+1}$  by  (i) computing $\rho = \rho^\text{NN} (\td_{j})$, (ii) random sampling of a bin index $m \, \in \{1, 2, ..., M \}$ in accordance with $\rho$, (iii) random sampling of $\tr_{j+1}$ from all $r^o_{i}$ in the $m$-th bin, $B_m$. 
Steps (ii) and (iii) can be combined: given  $\rho$, we sample $\tr_{j+1}$ randomly from all $r^o_{i}$, with weights $w_{i} = \sum_m  |B_m|^{-1} \, \rho_m \, \mathbf{1}(r^o_{i} \in B_m)$. Here $|B_m|$ denotes the number of training points $r^o_{i}$ in $B_m$, and $\mathbf{1}(.)$ is the indicator function. 

Training $\rho^\text{NN}$ can be seen as a task of probabilistic classification, with the bins $B_m$ as classes ($M$ in total). Focusing for now on the situation with $N=1$, a neural network for the classification task has $2J+2$ inputs (the dimension of the feature vector when $N=1$) and a $M$-dimensional softmax output layer. This is denoted a quantized softmax network (QSN). Note that this approach has been generalized to predict conditional kernel-density estimates for continuous variables, see \cite{ambrogioni2017kernel}.  However, here we stick to the discrete version.

\begin{figure}
    \centering
    \includegraphics[scale=0.5]{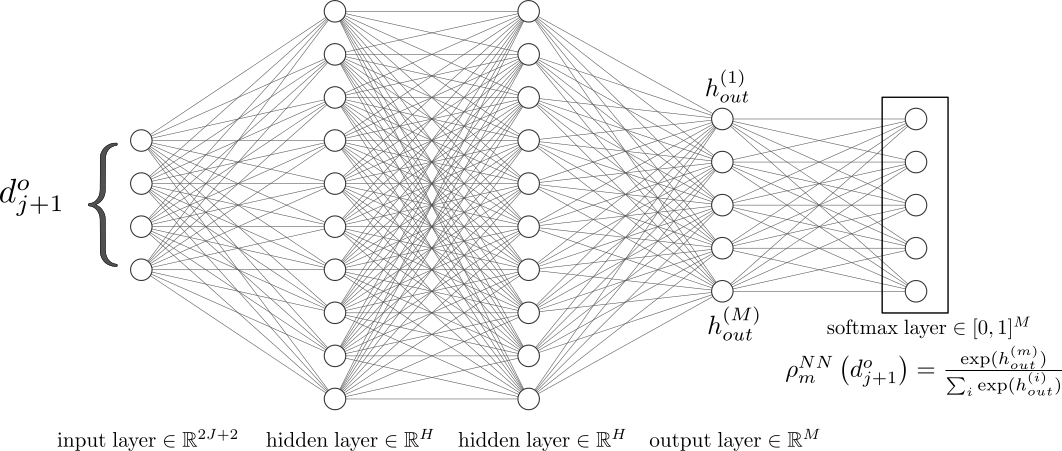}
    \caption{Diagram of a feed-forward QSN with two hidden layers. Here, $h^{(i)}_{out}$ is the output of the i-th neuron in the output layer, to which a softmax layer is attached.}
    \label{fig:nn}
\end{figure}

We use a feed-forward architecture (see Figure \ref{fig:nn}), noting that it can easily be replaced by other architectures in our set-up. The softmax layer computes the output probability mass function (pmf) as
\beq
\rho^\text{NN}_m (d_{j}^o) = \frac{\exp\left(h^{(m)}_{out}\right)}{\sum_{i=1}^M\exp\left(h_{out}^{(i)}\right)}.
\eeq
\noindent
Here, $h_{out}^{(i)}$ is the output of the i-th neuron in the output layer. Finally, a cross-entropy loss function \cite{aggarwal2018neural} is used for training the network:
\beq
\mathcal{L} = - \, \frac{1}{T-J-1} \, \sum_{j=J}^{T-1} \, \sum_{m=1}^M [\log \rho_m^\text{NN} (d_j^o)] \, \ind ( r^o_{j+1} \in B_m).
\label{eq:loss}
\eeq
\noindent
For a given $j$ index, the gradient of \eqref{eq:loss} with respect the output neurons (required to initialized a back propagation step), is given by $\partial \mathcal{L}/\partial h_{out}^{(m)} = \rho^{NN}_m - \ind ( r^o_{j+1} \in B_m)$ \cite{aggarwal2018neural}.

We emphasize that the method proposed here is fundamentally different from the approach of using a neural net to predict $\tr_{j+1}$ directly from $\td_j$ as in e.g. \cite{maulik2019subgrid, bolton2019applications,rasp2018deep}. In the latter approach, $\tr_{j+1}$ is modeled as a function of $\td_j$, with the neural net embodying the function as a deterministic mapping from the space of $\td_j$ to the space of $\tr_{j+1}$. By contrast, in our approach $\tr_{j+1}$ is resampled instead of modeled. Moreover, the mapping from the space of $\td_j$ to the space of $\tr_{j+1}$ in our approach is stochastic, not deterministic. Note that even if the neural net performs poorly, the $\tr_{j+1}$ that are generated remain consistent with the full model (as they are resampled from the full model), albeit possibly not well-matched with the feature vectors $\td_j$. Consistency in the case of vector-valued $r$ (i.e., $N>1$) is discussed next.

\subsection{Vector-valued $r$ and $x$}
\label{sec:vector}

To generalize from the case in which $r$ and $x$ are scalar ($N=1$) to the case where they are vector-valued ($N>1$), there are several possibilities. On the input side, the most straightforward approach is to enlarge the input layer from $2J+2$ to $2N(J+1)$ inputs. This simply reflects the increased dimension of the feature vector. 

\begin{figure}
    \centering
    \includegraphics[scale=0.6]{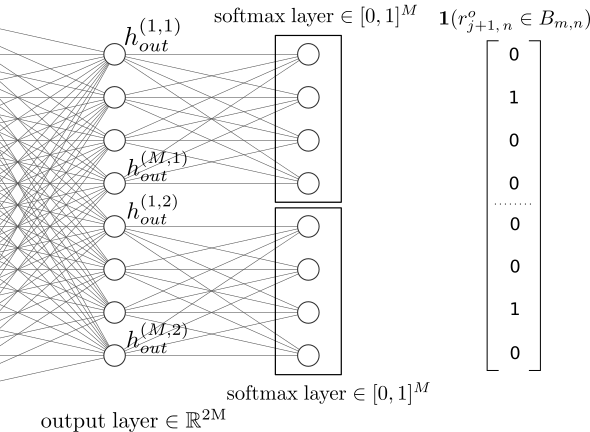}
    \caption{Diagram of the final layers of a feed-forward QSN with two probabilistic outputs with 4 bins each ($N=2$, $M=4$). Here, $h^{(m, n)}_{out}$ is the output of the m-th neuron in the output layer attached to the n-th softmax layer. Here, $\ind ( r^o_{j+1, n} \in B_{m, n})$ displays an example one-hot encoded data vector, for both softmax layers.}
    \label{fig:nn2}
\end{figure}

On the output side, if we bin the space in which the vector $r$ lives, we quickly run into curse of dimension problems again. Instead, we bin all individual elements of the vector $r$ separately. This scales linearly with dimension ($N$): with $N_{bins}$ bins for each vector element, we have $N \times N_{bins}$ bins in total. As a consequence, the output layer is no longer fully connected to a single softmax layer, see Figure \ref{fig:nn2} for a diagram involving $N=2$ probabilistic outputs. While each element of $r$ has its own independent softmax layer, all elements share the same feature vector and hidden layers. If there are dependencies between different elements of $r$ in the full model (\ref{eq:discrete}), they must be learned from the data $\{ r_j^o \}$. Such dependencies can be due to e.g. spatial correlations (with different vector elements representing different locations on a spatial grid). The separate elements of the vector $\tr_j$ generated this way are consistent with the full model, however their combination (i.e., the vector as a whole) may not be consistent if the network is poorly trained (in which case dependencies may be misrepresented).

While the scaling is linear with $N$, the number of output neurons can still become very large, especially for problems in two or three spatial dimensions (where $N$ can be e.g.\ $\mathcal{O}(10^6)$). In this case we can use other methods for discretizing the space of $r$. One such method is \emph{clustering} of the observations $\{r_j^o\}$, as also used in \cite{dorrestijn2013stochastic}. In this case, an entire vector $r_j^o$ is resampled at once from a selected cluster, rather than $N$ samples for the individual vector elements.  In this case the aforementioned consistency does not have to be learned from the data, but is guaranteed by the fact that all $r_j^o$ are produced by the full model (although again, the $\tilde{r}_{j+1}$ samples might still be ill-matched with the feature vectors $\td_j$).

Another method to tackle cases in which $N$ is set by the number of grid points on a spatial grid, is by formulating ``local" parameterizations. If $(x_j)_n$ denotes $x$ at time $t_j$ at grid point $n$ and similarly for $(r_j)_n$, a local parameterization to generate $(\tr_{j+1})_n$ takes only $(\tx_j)_n$ and $(\tr_{j})_n$ (and possibly their histories) as its inputs. Thus, it ignores $(\tx_j)_{n'}$ and $(\tr_{j})_{n'}$ at other grid points $n' \neq n$. This is standard practice for parameterizations in e.g. atmosphere-ocean science. However, it may fail to capture spatial correlations.

Returning to our approach with $N$ softmax layers,
the loss function is also slightly modified, as it now includes summation over the $N$ layers: 
\beq
\mathcal{L} = - \, \frac{1}{T - J - 1} \, \sum_{j=J}^{T-1} \, \sum_{n=1}^N \, \sum_{m=1}^M [\log \rho_{m, n}^\text{NN} (d_j^o)] \, \ind ( r^o_{j+1,\,n} \in B_{m,\,n}). 
\label{eq:loss2}
\eeq
\noindent
Here, $\rho^{NN}_{m, n}$ is the QSN predicted probability for the m-th bin of the n-th softmax layer, and the indicator function $\ind ( r^o_{j+1, n} \in B_{m, n})$ represents the one-hot encoded data per softmax layer (see Figure \ref{fig:nn2}). The gradient of \eqref{eq:loss2} with respect to the output neurons retains the same expression as in the case of a single softmax layer, i.e.\ $\partial \mathcal{L}/\partial h_{out}^{(m, n)} = \rho^{NN}_{m,n} - \ind ( r^o_{j+1,\, n} \in B_{m\,,n})$.

\section{Numerical experiments}


In this section we present several examples to test the approach proposed in the previous section: model reduction according to \eqref{eq:reduced2} with the resampling step in \eqref{eq:reduced2b} implemented using  QSNs. In these tests, it is not our aim to recreate the exact trajectories of the original coupled model \eqref{eq:discrete} with the reduced model. The loss of information due to model reduction, the stochastic nature of our surrogates, and the nonlinearity of multiscale problems, makes this impossible after a certain time integration length. Instead, our goal for the setup \eqref{eq:reduced2} is to reproduce the time-averaged statistical properties of the original macroscopic variables $x$, e.g,.\ the probability density function or the auto-correlation function.

\subsection{Model equations Lorenz 96 \label{sec:l96}}

\begin{figure}
    \centering
    \includegraphics[scale=0.6]{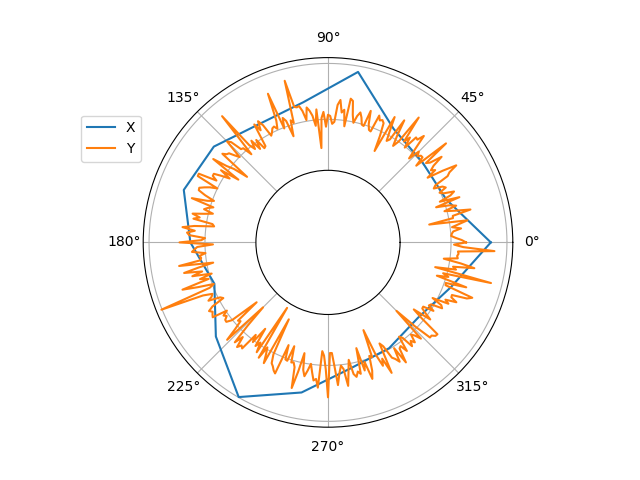}
    \caption{A snapshot of the solution of the two-layer L96 system \eqref{eq:lorenz96}. Here $X$ denotes the vector of all $X_n$ solutions and $Y$ is the vector containing all $Y_{l,n}$.}
    \label{fig:l96}
\end{figure}

As a proof of concept, we test our setup on the well-known two-layer Lorenz 96 (L96) system, originally proposed by \cite{lorenz1996predictability}, which is a toy model for the atmosphere. It consists of a set of $N$ ODEs describing the evolution of the macroscopic variables $X_n$, of which each ODE is coupled to $L$ microscopic variables $Y_{l, n}$ (note that $n$ and $l$ are spatial grid indices here):
\begin{align}
    \frac{\dif X_n}{\dif t} = 
    X_{n-1}\left(X_{n+1} - X_{n-2}\right) - X_n - F + r_n \nonumber\\
    r_n := \frac{h_x}{L}\sum_{l=1}^L Y_{l,n} \nonumber \\
    \frac{\dif Y_{l, n}}{\dif t} = 
    \frac{1}{\epsilon}\left[Y_{l+1, n}\left(Y_{l-1, n} - Y_{l+2, n}\right) - Y_{l, n} + h_yX_n\right].
    \label{eq:lorenz96}
\end{align}
The macroscopic and microscopic variables $X_n$ and $Y_{l,n}$ are considered variables on a circle of constant latitude (see Figure \ref{fig:l96}), where the indices $n=1,\cdots,N$ and $l=1,\cdots,L$ denote the spatial location. Note that there are $L$ microscopic solutions $Y_{l,n}$ for every $n$, such that $Y=(Y_{1,1}, Y_{1,2},\cdots, Y_{L,N})\in\mathbb{R}^K$, where $K=NL$. The circular shape of the domain is imposed through periodic boundary conditions:
\begin{align}
 X_{n} = X_{n+N}, \quad Y_{l, n} = Y_{l,\;n+N}, \quad Y_{l+L,\;n} = Y_{l, \;n+1}.
\end{align}
\noindent
The two-step Adams–Bashforth method is chosen as the numerical time-discretization procedure, and we set $\Delta t = 0.01$ as the time step. Note that \eqref{eq:lorenz96} is only used to generate the training data. In line with \eqref{eq:reduced2}, we only solve the $X_n$ equation with a QSN surrogate for $r_n$ when predicting.

%

\begin{figure}
    \centering
    \includegraphics[scale=0.4]{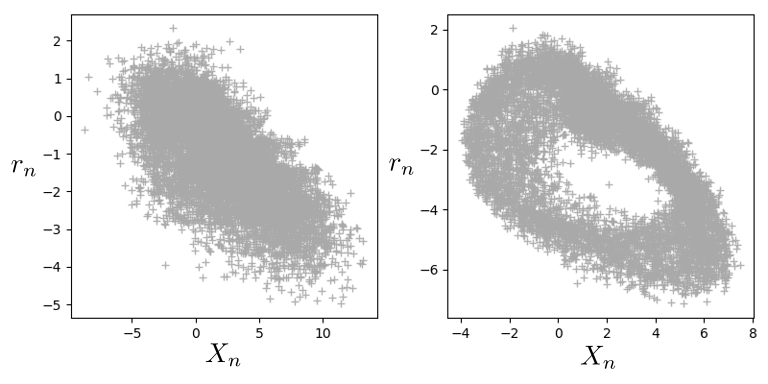}
    \caption{Scatter plot of $X_n$ vs $r_n$ for $\{N, L, F, h_x, h_y, \epsilon\} = \{18, 20, 10, -1, 1, 0.5\}$ (left), and $\{N, L, F, h_x, h_y, \epsilon\} = \{18, 20, 10, -2, 1, 0.5\}$ (right).}
    \label{fig:uni_vs_bi}
\end{figure}

The behavior of the system is governed by the parameter settings of $\{N, L, F, h_x, h_y, \epsilon\}$. Commonly, $\epsilon$ is chosen such that $\epsilon\ll1$, in which case a clear time-scale separation between the macroscopic and microscopic variables is imposed. We will follow \cite{crommelin2008subgrid} and use $\epsilon=0.5$, such that no clear temporal gap exists, which is more realistic for turbulent geophysical flows, and more challenging for parameterizations. Specifically, we will use two different parameter settings:
\begin{enumerate}
    \item The unimodal setting: $\{N, L, F, h_x, h_y, \epsilon\} = \{18, 20, 10, -1, 1, 0.5\}$,
    \item The bimodal setting: $\{N, L, F, h_x, h_y, \epsilon\} = \{18, 20, 10, -2, 1, 0.5\}$.
\end{enumerate}
\noindent
The naming convention stems from the number of modes in the probability density functions (pdfs) of the $X_k$ variables. The bimodal setting intensifies the feedback from the microscopic ODEs to the macroscopic ODE, by modifying the $h_x$ parameter (see \eqref{eq:lorenz96}). Specifically, we decrease $h_x$ from -1 to -2. In the unimodal setting, both the $X_n$ and $r_n$ pdfs are unimodal and nearly Gaussian, see Figure \ref{fig:l96_pdf}. This is no longer the case for $h_x = -2$, when the pdfs become non-symmetric and bimodal. The difference between the two parameter settings can also clearly be seen from a scatter plot of $X_n$ vs $r_n$, see Figure \ref{fig:uni_vs_bi}.

Other authors have used L96 as a benchmark for machine learning methods. For instance in \cite{chattopadhyay2019data}, three separate methods are used to predict the full $X^{(j+1)}=(X_1(t_{j+1}),\cdots, X_N(t_{j+1}))$ vector. Also, the authors of \cite{dueben2018challenges} used a neural network to predict the tendency $\Delta X = X^{(j+1)} - X^{(j)}$ one time step ahead. Note that these are deterministic approaches, and are not applied in the context of parameterization (they are used to predict $X^{(j+1)}$ itself rather than $r^{(j+1)}$). Closer to our approach, the interesting recent work of \cite{christensen2020machine} uses conditional Generative Adversarial Networks (GANs) for the purpose of stochastic parameterization of the L96 subgrid-scale term. GANs have a different architecture than our QSNs, and unlike our approach, include a stochastic component in the inputs. The approach from \cite{christensen2020machine} does not involve resampling, and was tested on the L96 model at a (single) different parameter setting. 

\subsection{Learning procedure \label{sec:learning}}

To inform the weights of a QSN, we use back propagation with Stochastic Gradient Descent with the RMSProp optimizer and a learning rate of 0.001 \cite{aggarwal2018neural}. The number of training iterations was set to 10000. After experimenting with different networks, we selected 3 hidden layers, each with 256 neurons and leaky Rectified Linear Unit activation functions. The output layer, which feeds into the softmax layers, is linear. Furthermore, the input features $d^o_j$ are standardized to have zero mean and a standard deviation of one, before being fed into the QSN. We will create both a local surrogate (trained on a single spatial location), and a stochastic surrogate for the full vector-valued $r:=[r_1,\cdots, r_N]^T$, in which case we have $N$ softmax layers, see Section \ref{sec:vector}.

A common practice in machine learning is to leave out a part of the data set (i.e.\ to not use it in training), to test the accuracy of the final model.  In our case however, the trained neural network itself is not the final model, it is merely a source term in the final model (the macroscopic ODE \eqref{eq:reduced2a}). To test the accuracy of the neural network, we therefore have to perform a simulation with two-way coupling between the ODE and the neural net. We leave out the final 50 \% of the data and test the ability of this coupled system to predict the macroscopic statistics of the test set.

\subsection{Results: verification}

\begin{figure}
    \centering
    \includegraphics[scale=0.7]{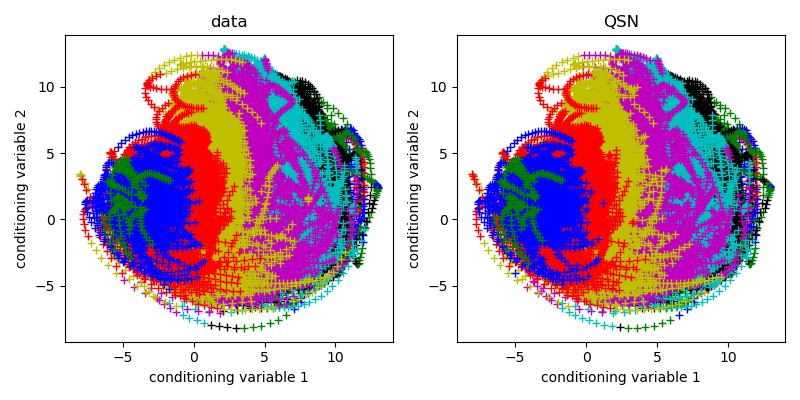}
    \caption{A twice time-lagged $X_n$ conditioning feature is shown on both axes, and the color of the ``+'' symbols denotes the index of the bin with highest probability of the corresponding $r_{j+1}$ sample. The left plot shows the output bin indices of the training data, and the right plot show the predicted indices corresponding to QSN output pmf, when using the training features $d^o_j$ as input. 
    }
    \label{fig:bins_L96}
\end{figure}

To visualize the complexity of the bin classification and to verify quality of a trained QSN, consider Figure \ref{fig:bins_L96}. Here we show a scatter plot with a lagged $X_n$ feature on both axes. The symbols in the scatter plot are color coded for the corresponding $r$ bin index. We show both the exact results from the training data and the QSN prediction, using the training $d^o_j$ as input features. Figure \ref{fig:bins_L96} verifies that the QSN can learn a good bin index representation in the space of lagged input variables. A typical misclassification error for each softmax layer $n=1,\cdots, N$ (defined as $\mathrm{arg}_m\mathrm{max} \;\rho^{NN}_{m, n}(d^o_j) \neq \mathrm{arg}_m\mathrm{max}\; \ind ( r^o_{j+1, n} \in B_{m,n})$), is roughly 3-4\%.

\begin{figure}
    \centering
    \includegraphics[scale=0.37]{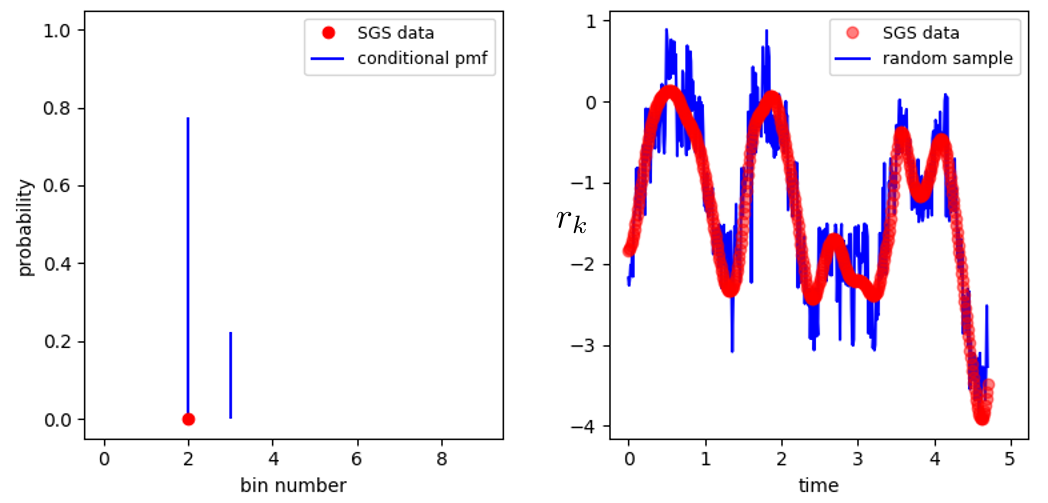}
    \caption{At any given time $t_j$, the QSN predicts a conditional pmf (left). The predicted bin index from this pmf feeds into our data resampling software, which outputs a random $\tilde{r}_{j+1}$. Its time series (for a single spatial point $n$), along with the subgrid-scale (SGS) data $r$ are shown on the right. \label{fig:qsn}}
\end{figure}

To verify the random resampling software we plot the time series of the data $r^o_{j+1}$ and stochastic estimate $\tilde{r}_{j+1}(d^o_j)$ in Figure \ref{fig:qsn}.

\subsection{Results: validation}

We simulate the reduced system \eqref{eq:reduced2} from $t = 0$ until $t = 1000$, while the training data spanned $t\in[0, 500]$. Our macroscopic statistics of interests are the probability density function and auto-correlation function of $X_n$, as well as the cross-correlation function of $X_n$ and $X_{n+1}$. In addition, we compute the same statistics for $r_n$. All statistics are averaged over all $n$, and computed on the test set only, i.e. using data from $t\in[500, 1000]$. 

Our main goal here is twofold: We investigate the importance of the length of the history (memory depth) in the feature vector ($J$ in definition (\eqref{eq:features})). Furthermore, we demonstrate the relevance of a stochastic approach by comparing resampling with using averages (as discussed in section \ref{sec:stoch}). In addition, we compare the performance of a surrogate that is trained and applied locally, and a surrogate which predicts the entire $\tilde{r}_{j+1}\in\mathbb{R}^N$ vector at once. 

\subsubsection{Short vs long memory}

\begin{figure}
    \centering
    {\includegraphics[scale=0.5]{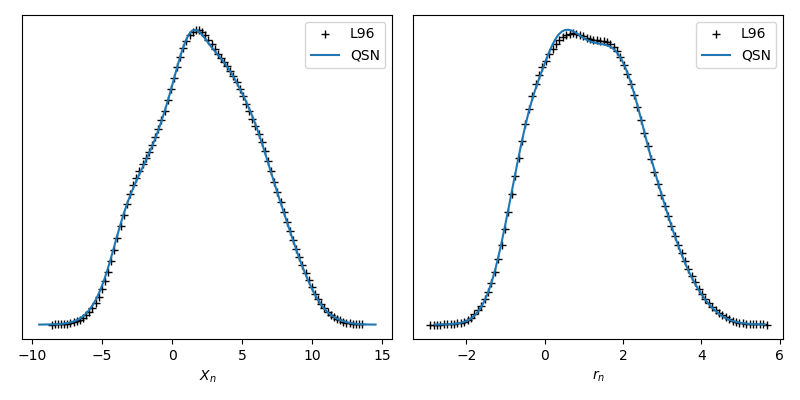}}
    \caption{Probability density functions of $X_n$ (left) and $r_n$ (right). Dots denote the full two-layer solution and solid lines denote the reduced model with QSN surrogate. \label{fig:l96_pdf}}
\end{figure}

\begin{figure}  
    \centering
    {\includegraphics[scale=0.5]{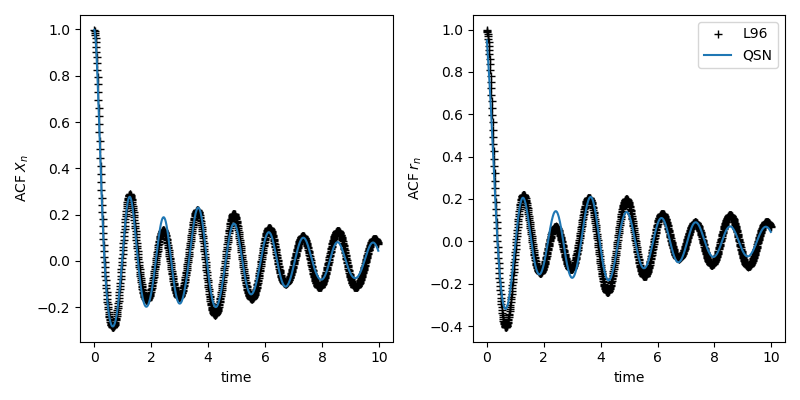}}
    \caption{Auto-correlation function of $X_n$ (left) and $r_n$ (right). Dots denote the full two-layer solution and solid lines denote the reduced model with QSN surrogate. \label{fig:l96_acf}}
\end{figure}

\begin{figure} 
    \centering
    {\includegraphics[scale=0.5]{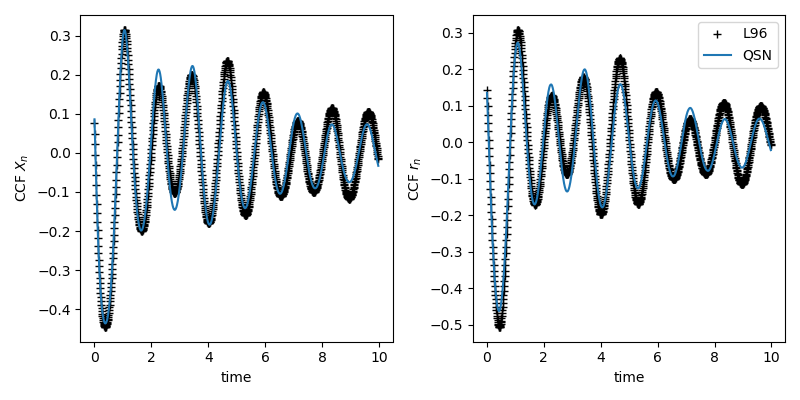}}
    \caption{Cross-correlation function of $X_n$ and $X_{n+1}$ (left) and $r_n$ and $r_{n+1}$ (right). Dots denote the full two-layer solution and solid lines denote the reduced model with QSN surrogate. \label{fig:l96_ccf}}
\end{figure}

For our first test case we will the unimodal parameter values. We first create a full-vector QSN with a short memory, e.g.\ with two time-lagged $X:=(X_1, \cdots, X_N)$ vectors in $d^o_j$. Let $X^{(j)}:=X(t_j)$. The statistical results, when using $d_j^o=(X^{(j)}, X^{(j-9)})$, are shown in Figures \ref{fig:l96_pdf}-\ref{fig:l96_ccf}. Despite the short memory in $d^o_j$, these display a good match between the statistics of the reduced and the full L96 system. Note that in this case $d_j^o\in\mathbb{R}^{2N}$, which would be impossible in the case of a surrogate that bins the input space, as in \cite{verheul2016data}. 

\begin{figure}
    \centering
    {\includegraphics[scale=0.5]{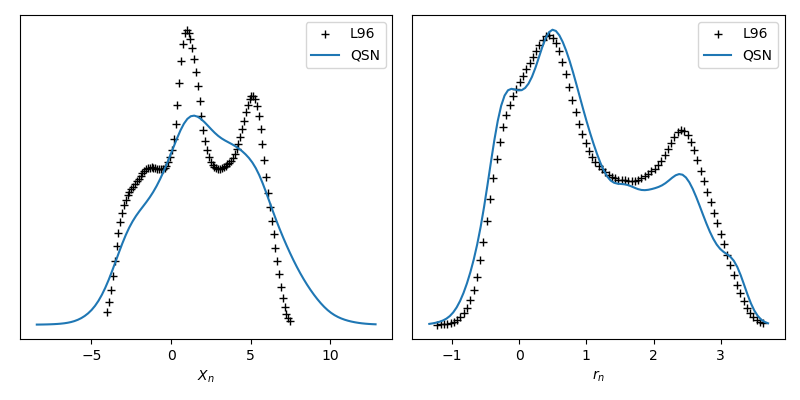}}
    \caption{Probability density functions of $X_n$ (left) and $r_n$ (right), using $d^o_j = (X^{(j)}, X^{(j-1)}, \cdots, X^{(j-9)})$. Plus symbols denote the full two-layer solution and solid lines denote the reduced model with QSN surrogate. \label{fig:l96_pdf_lag10}}
\end{figure}
\begin{figure}  
    \centering
    {\includegraphics[scale=0.5]{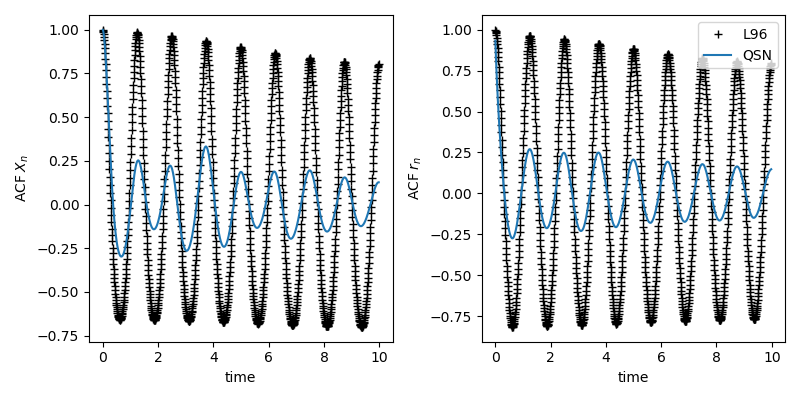}}
    \caption{Auto-correlation function of $X_n$ (left) and $r_n$ (right), using $d^o_j = (X^{(j)}, X^{(j-1)}, \cdots, X^{(j-9)})$. Plus symbols denote the full two-layer solution and solid lines denote the reduced model with QSN surrogate. \label{fig:l96_acf_lag10}}
\end{figure}
\begin{figure} 
    \centering
    {\includegraphics[scale=0.5]{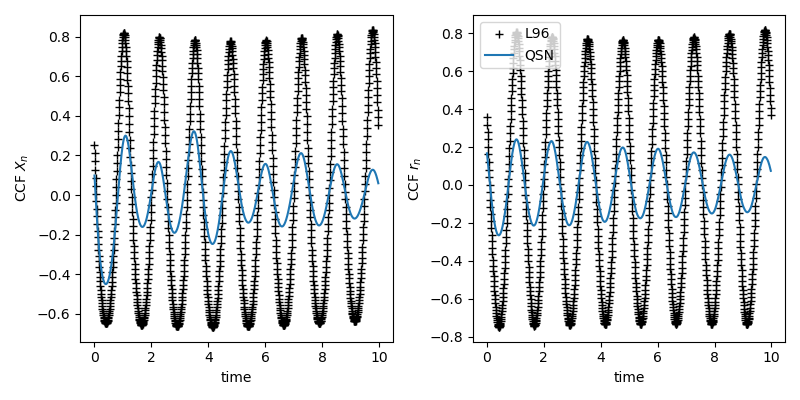}}
    \caption{Cross-correlation function of $X_n$ and $X_{n+1}$ (left) and $r_n$ and $r_{n+1}$ (right), using $d^o_j = (X^{(j)}, X^{(j-1)}, \cdots, X^{(j-9)})$. Plus symbols denote the full two-layer solution and solid lines denote the reduced model with QSN surrogate. \label{fig:l96_ccf_lag10}}
\end{figure}
\begin{figure}
    \centering
    {\includegraphics[scale=0.5]{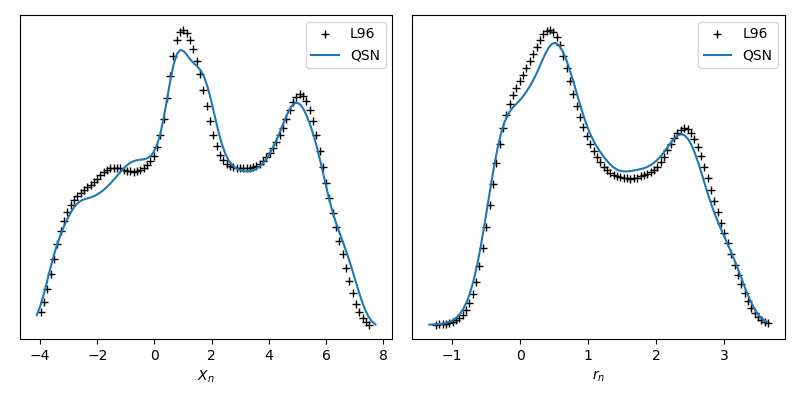}}
    \caption{Probability density functions of $X_n$ (left) and $r_n$ (right), using $d^o_j = (X^{(j)}, X^{(j-1)}, \cdots, X^{(j-74)})$. Plus symbols denote the full two-layer solution and solid lines denote the reduced model with QSN surrogate. \label{fig:l96_pdf_lag75}}
\end{figure}
\begin{figure}  
    \centering
    {\includegraphics[scale=0.5]{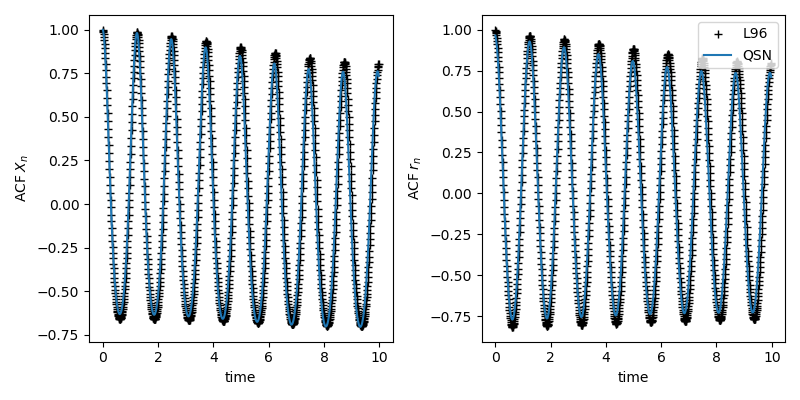}}
    \caption{Auto-correlation function of $X_n$ (left) and $r_n$ (right), using $d^o_j = (X^{(j)}, X^{(j-1)}, \cdots, X^{(j-74)})$. Plus symbols denote the full two-layer solution and solid lines denote the reduced model with QSN surrogate. \label{fig:l96_acf_lag75}}
\end{figure}
\begin{figure} 
    \centering
    {\includegraphics[scale=0.5]{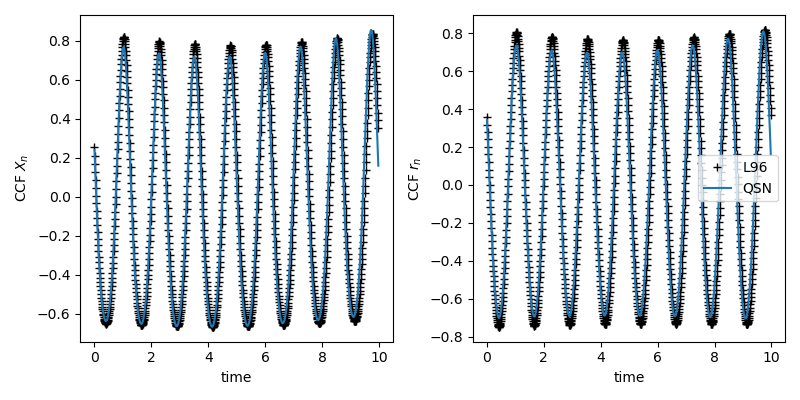}}
    \caption{Cross-correlation function of $X_n$ and $X_{n+1}$ (left) and $r_k$ and $r_{k+1}$ (right), using $d^o_j = (X^{(j)}, X^{(j-1)}, \cdots, X^{(j-74)})$. Plus symbols denote the full two-layer solution and solid lines denote the reduced model with QSN surrogate. \label{fig:l96_ccf_lag75}}
\end{figure}

Above, we used a parameter setting for the L96 model that makes parameterization challenging because of the lack of scale separation ($\epsilon=0.5$ as discussed). 
Here, we put more strain on the QSN by employing the bimodal parameter setting. This is a harder test case, and more care needs to be taken with the specification of the feature vector $d^o_j$. In fact, the amount of memory in $d^o_j$ becomes very important. If we use 10 $X$ vectors, i.e. $d^o_j = (X^{(j)}, X^{(j-1)}, \cdots, X^{(j-9)})$, we obtain the results of Figures \ref{fig:l96_pdf_lag10}-\ref{fig:l96_ccf_lag10}, which display a clear failure of capturing the reference statistics. Also note the bimodal nature of the reference pdfs. We performed further tests with 25, 50 and 75 lagged $X$ vectors, and only obtained good statistics with 75 vectors. These results are shown in Figures \ref{fig:l96_pdf_lag75}-\ref{fig:l96_ccf_lag75}.





\subsubsection{Stochastic vs deterministic}

We recall \eqref{eq:reduced2b} and the discussion in section \ref{sec:stoch} about resampling versus using bin averages. Here we compare sampling from the set $S_{j+1}$ (i.e.\ the output bin in the case of a QSN), with using the sample mean of $S_{j+1}$. If we also use $\mathrm{\arg}_m\mathrm{max}\;\rho^{NN}(\tilde{d}_j)$ as the predicted output bin index, we obtain a completely deterministic surrogate $\tilde{r}_{j+1}$. We tested this approach on the full-vector surrogates of the preceding section, and obtained similar results as with resampling. However, when applying the surrogate locally, we find a significant impact due to the stochastic nature of $\tilde{r}_{j+1}$. As mentioned, the local surrogate is trained on a single spatial location, and during prediction it is applied independently for each $X_n$ equation. Such an approach is not uncommon, see e.g. \cite{crommelin2008subgrid, christensen2020machine, rasp2019online}, and it matches the local nature of traditional parameterization schemes.

\begin{figure}  
    \centering
    {\includegraphics[scale=0.5]{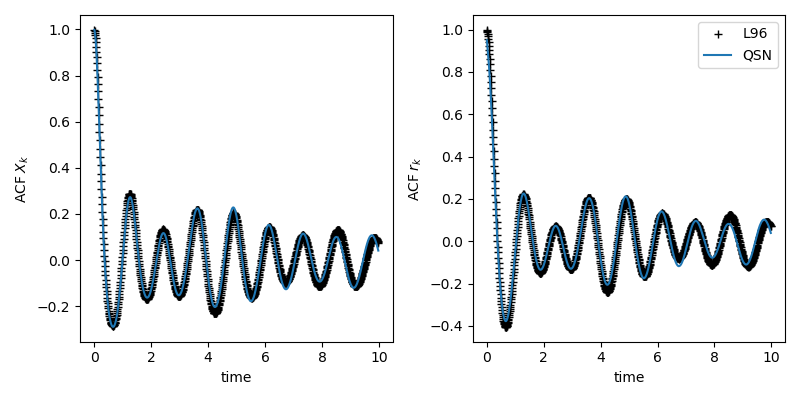}}
    \caption{Auto-correlation function of $X_n$ (left) and $r_n$ (right), using $d^o_j = (X^{(j)}, X^{(j-1)}, \cdots, X^{(j-74)})$. Plus symbols denote the full two-layer solution and solid lines denote the reduced model with the {\it local and stochastic} QSN surrogate. \label{fig:l96_acf_local}}
\end{figure}

\begin{figure}  
    \centering
    {\includegraphics[scale=0.5]{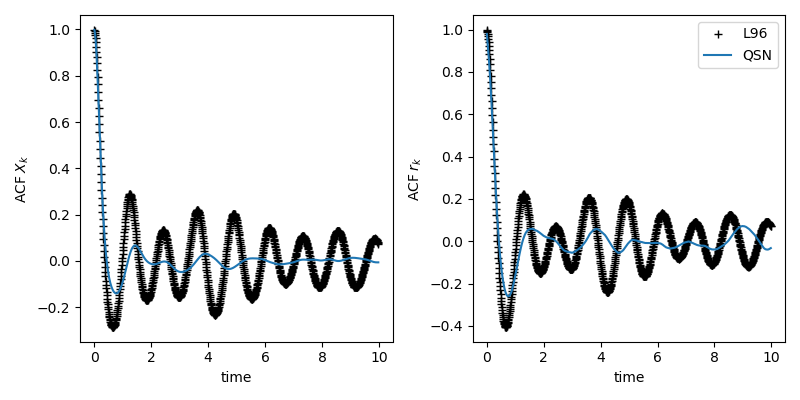}}
    \caption{Auto-correlation function of $X_n$ (left) and $r_n$ (right), using $d^o_j = (X^{(j)}, X^{(j-1)}, \cdots, X^{(j-74)})$. Plus symbols denote the full two-layer solution and solid lines denote the reduced model with the {\it local and deterministic} QSN surrogate. \label{fig:l96_acf_local_det}}
\end{figure}

As an example, consider the auto-correlation results of Figures \ref{fig:l96_acf_local}-\ref{fig:l96_acf_local_det}, which show the results of the stochastic and deterministic surrogate respectively. Note that here, $h_x = -1$ was used, i.e.\ the easier of the two parameter setting we consider. The stochastic surrogate clearly outperforms its deterministic counterpart in this case. Other statistics (acf, pdf) showed similar results.


\subsubsection{$X$-only vs $X-r$ conditioning}

In Section \ref{sec:stoch} we outlined the general case in which the feature vector consists of time-lagged $x$ and $r$ variables. Yet, thus far we have only shown $X$-only surrogates for the L96 system. We found that including $r_j$ in $d^o_j$ can contribute to obtain a small training error, especially in the case of local surrogates. However, we did not obtain robust statistical predictions when doing so. These results are in line with those from \cite{christensen2020machine}. As mentioned in Section \ref{sec:l96}, these authors considered stochastic parameterizations for L96 by means of conditional GANs. For a variety of tests, their X-only GANs clearly outperformed the GANs conditioned on both $X$ and $r$. One possible cause may be overfitting, due to the very strong correlation between $r_j$ and $r_{j+1}$ \cite{christensen2020machine}.

\subsection{Software}

The source code used to generate the results can be downloaded from \cite{edeling2020physd}.

\section{Discussion and future challenges}

The preceding results show that, for the L96 system, the amount of memory ($J$) in $d^o_j$ was not of major importance in the tests with the unimodal parameter settings ($h_x=-1$). However, the more challenging bimodal parameter settings ($h_x = -2$, right subplot of Figure \ref{fig:uni_vs_bi}) resulted in a problem for which memory became crucial. In general, we expect that for more complicated (geophysical) flow problems, memory will play an important part. It is clear however, that the ``optimal" $d^o_j$ is problem dependent, and a systematic procedure for designing the best feature vector is an interesting avenue for future research. This could involve changing the network architecture (e.g.\ combining resampling with Long Short-Term Memory networks \cite{hochreiter1997long}), or finding optimal time lags using approaches as described in \cite{frank2001time}.


Another clear avenue for future research is to apply machine learning with resampling as proposed here to more complex flows. In Section \ref{sec:vector}, we have discussed ways to deal with the large output dimension, including clustering of observations $\{r_j^o\}$, as in \cite{dorrestijn2013stochastic}. An interesting test problem is a two-dimensional ocean model, as in e.g.\ \cite{berloff2005random, verheul2017covariate,edeling2019towards}.

Finally, as mentioned in Section \ref{sec:learning}, we train the QSN separately on the data, and afterwards a validation procedure involves a two-way coupling between the ODEs and the QSN. This gave satisfactory results for the L96 model, but for more complicated problems, such an ``offline" training strategy could lead to instabilities in the ``online", two-way coupled simulation as discussed in \cite{rasp2019online}. Developing new learning procedures, in which the neural network is trained while it is part of the larger dynamical system is of interest.

\section{Conclusion \label{sec:con}}

We presented a machine-learning method for the conditional resampling of subgrid-scale data of multiscale dynamical systems, resulting in a stochastic parameterization for the unresolved scales. The current model is comprised of a feed-forward architecture with (multiple) softmax layers attached to the output. The output data is divided into a finite number of non-overlapping intervals (denoted as `bins'), and the softmax layers predict a discrete probability density function over these bins, conditioned on time-lagged macroscopic input features. First a bin is sampled from this distribution, which is followed by randomly selecting a reference subgrid-scale data point from the identified bin. This stochastic surrogate model then replaces the original subgrid-scale term in the dynamical system, and we validate the method by examining the ability of this system to capture the long-term statistics of the resolved, macroscopic variables.

In this initial study we considered the Lorenz 96 system at two different parameter settings. We were able to validate our method on this problem, provided that a sufficient number of time-lagged variables were included in the feature vector. Finally, we also found that overall, the stochastic nature of the surrogate led to more robust performance.

\section*{Acknowledgements}
This research is funded by the Netherlands Organization for Scientific Research (NWO) through the Vidi project "Stochastic
models for unresolved scales in geophysical flows", and from the European Union Horizon 2020 research and innovation programme under grant agreement \#800925 (VECMA project). 

\bibliographystyle{plain}
\bibliography{sample.bib}

\end{document}